# STRONG REPETITIONS IN STURMIAN STRINGS

Ayşe Karaman
makalelr@yahoo.ca

**ABSTRACT:** In this paper, we analyze the periodic factors of Sturmian words for the findings to lead to a linear-time algorithm for the computation of runs in this class of words which, to our best knowledge, is an open problem in literature.

## 1. INTRODUCTION

A *Sturmian string*, or a *Sturmian word* is a binary word infinite at the right end which has the minimum number of distinct factors and which is not ultimately periodic. Properties and representations of Sturmian words have been studied by various authors including [AS03, B96, BS93, CH73, DM94, K98, L02]. A *strong repetition*, or shortly a *repetition* is the concatenation of a finite word with itself at least twice. A *run* is a maximal repetition in a word so that it is not a proper factor of another repetition in the same word. [FKS00] studied the computation of repetitions in Sturmian words. Prior to that, [IMS97] presented a linear-time algorithm to compute the runs in Fibonacci words. [DL03] studied the repetitions for their length and powers.

In this paper, we characterize the runs in Sturmian words and show their computation in linear time. Our analysis uses a formulation of Sturmian words introduced in [K98] as a powerful to study the combinatorial structures of this class of words. Following the preliminaries and a literature survey in the next section, we present our results in Section 3. The last section concludes with additional remarks on the computations based on the findings of this paper.

## 2. PRELIMINARIES

A Sturmian word is a binary word infinite at one end which has minimum number of distinct factors and which is not ultimately periodic. In this text, we use $A=\{a,b\}$ for the binary alphabet. Another way to identify a Sturmian word is to see whether it is balanced while not ultimately periodic. A binary word is defined as balanced if for every finite factor pair $u$ and $v$ of equal length, the number of letter $b$'s (equivalently $a$'s) in $u$ and $v$ differ by at most 1.

Consider the morphism $\alpha_i(a)=a^{p_i}b$, $\alpha_i(b)=a^{p'_i}b$ where $p_i, p'_i \geq 1$, $|p_i - p'_i|=1$. This morphism is clearly Sturmian so that a word is balanced and not ultimately periodic if and only if its image under such a morphism is (see for example [B96]). In [K98], we used this class of $\alpha$-morphisms to characterize the Sturmian words. Let $\pi=\{(p_1, p'_1), (p_2, p'_2), ...\}$ be a *defining sequence* with each $(p_i, p'_i)$ a *valid parameter pair* satisfying $p_i, p'_i \geq 1$, $|p_i - p'_i|=1$. Let $\alpha_\pi(a)$ denote

$$\alpha_\pi = ... \bullet \alpha_i ... \bullet \alpha_2 \bullet \alpha_1.$$

[K98] shows that a finite or infinite Sturmian word is a factor of a word $\alpha_\pi(a)$ for some defining sequence $\pi$, and every $\alpha_\pi(a)$ generates a Sturmian word so that the class of Sturmian words can be formulated by defining sequences in the $\alpha$-morphic representation. This formulation is a powerful tool to analyze the structure of Sturmian words as composition of levels of identical properties through which the patterns can be detected and their iteration can be computed at the ultimate level. In this paper, we use this tool to analyze the runs in Sturmian words for their characterization.

Consider a Sturmian word $S_i = \alpha_\pi(a) = \alpha_i(\alpha_{i-1}( \ldots (\alpha_1(a)))$ on a finite defining sequence $\pi$. We say that $S_i$ is of *level i*. Observe that a Sturmian word $S_i$ of level $i$ is a concatenation of blocks $a^{p_i} b$ and $a^{p'_i} b$. We call each such block, namely $\alpha_i(a)$ and $\alpha_i(b)$ *full blocks*. A factor in $S_i$ which is a proper factor of a full block is called a *partial block*. In $S_i$, let $p_{min} = \min\{ p_i, p'_i \}$ and correspondingly $p_{max} = p_{min}+1$. We name the full blocks $a^{p_{min}} b$ and $a^{p_{max}} b$ respectively as *short* and *long blocks*. Observe in such a case that $a^{p_{min}} b$ is a full when it follows a full block occurrence, a partial block otherwise when it is a proper suffix of a long block. In $S_i$, $\alpha_i(a)$ is the *repeating block*. A factor of $S_i$ which is a concatenation of full blocks is *block-complete*. We define $\alpha_i^{-1}(s)$ over a defining sequence $\pi$ to denote the inverse morphism on a block-complete word $s$ to map full blocks into words so that $\alpha_i(\alpha_i^{-1}(s)) = s$ (see [K98] for a detailed definition of $\alpha_i^{-1}$). We use the notation $s[i]$ and $s[i..j]$, $i \leq j$ to denote respectively the letter in position $i$ and the factor in positions $i$ through $j$ of word $s$. $|s|$ denotes the length of a word $s$. Given a finite word $u$, we name a word $u[i+1..|u|]u[1..i]$ for any value of $i=1..|u|-1$ as a *circular shift* of $u$. Let $u=w[i..i+|u|-1]$ be a factor of a word $w$. Then, $(n)$-*left shift* of $u$ in $w$ is defined as $w[i-n..i+|u|-1-n]$. Similarly, $(n)$-*right shift* of $u$ in $w$ refers to $w[i+n..i+|u|-1+n]$. When we use the term *left* (or *right*) shift of $u$, we refer to an $(n)$-left (or -right) shift of $u$ in the word for an arbitrary value of $n=1..|u|-1$.

A *strong repetition* or a *repetition* is a word $u=g^k$ where the period $g$ is a finite non-empty word over the alphabet, and $k \geq 1$, $k$ is an integer. A repetition is *primitive* if its period is primitive, i.e. is not a repetition itself. In this text, repetition refers to a primitive repetition unless otherwise stated. For a repetition $u=g^k$ which is a subword of $S_i$ so that $S_i = cud$ and $c[|c|-|g|+1..|c|] \neq g$, the longest common non-empty suffix of $c$ and $u$ is called the *left extension* of $u$. Similarly, the longest common non-empty prefix of $d$ and $u$ where this prefix is shorter in length than $g$ is called the *right extension* of $u$. We use the notation $(g, k, e_1, e_2)$ to denote a maximal repetition, or a *run* $g^k$ with respectively the left and right extensions $e_1$ and $e_2$. Observe that a run $(g, k, e_1, e_2)$ satisfying $|e_1|+|e_2|<|g|$ can be normalized into the form of a triple $(g', k', e')$ where $g'=g[|g|-|e_1|+1..|g|]g[1..|g|-|e_1|]$, $e'=g'[1.. |e_1 e_2|]$, $k'=k$. The *index* of a word $w$ is a fractional value which is the maximum of the powers of runs in a word $w$. We refer by *integer index* to maximum of the integer powers of runs in the word.

Repetitive sequences in words are researched from a number of viewpoints. See Chapter 1 of [AS03] for the problem descriptions on repetitive sequences. See [AS03, CIR09] for literature survey of this topic under the general field of Word Combinatorics. Following are the results in literature relating to the characterization of runs in Sturmian words, which is the context of this study:

- *Characterization of Runs*: Iliopoulos and Smyth [IMS97] computed the runs in the standard Fibonacci word for their periods, positions and powers in the word by using a linear-time algorithm. Franĕk et. al [FKS00] studied the runs for their positions in a given Sturmian word and presented an algorithm to compute the repetitive sequences at such positions. Damanik and Lenz [DL03] characterized the integer powers in the Sturmian words for their existance for each period length in the word and, for each integer power $k$ admitting a repetition of the form $w^k$, enumerated the periods $w$. Note that in [DL03], the problem domain is not restricted to the maximal and primitive repetitions as in runs. The repetitive factors of runs are analyzed across their varying periods, circular shifts of periods and period lengths, and each run is $(p,t,e)$ is enumerated as $\theta(|p|t*ln(t))$ distinct repetitions without an explicit specification of runs in the word. Cao and Wen [CW03] independently showed some of the results in [DL03] in their comprehensive study on the singular words. Specifically, [CW03] characterized the repetitions of periods which are circular shifts of the finite standard words. They also showed that no period of length unequal to the length of a standard word has a singular factor. De Luca and De Luca [DD06] showed that the finite Sturmian factors are exactly the set of words whose minimum fractional periods are circular shifts of standard words.
- *Number of Runs*: [DL03] computed the number of squares of distinct periods and the total number of occurrences of squares in a Sturmian word given by its slope, which extended a result studied on Fibonacci words by Fraenkel and Simpson [FS99]. [PR09] use the results in [DL03] to study the asymptotic behaviour of squares in Sturmian words. Baturo et. al [BPR08] developed an explicit formula for the number occurrences of runs in a finite standard word given by its continued fraction expansion and showed that this count has an upper bound of 4/5 of the word length. The results in [BPR08] also verify the results in [CI07, KK99, PSS08, R06] (cf. [CIR09]) for the class of Sturmian words that the number of runs in an arbitrary word is linear in the length of the word. In relevance to the topic, de Luca and De Luca [DD06] give a formula to return the number of primitive Sturmian words of a given length.
- *Index of Sturmian words*: [K83] showed that the integer index of the Fibonacci word is 3. Mignosi [M89] showed that a Sturmian word has a finite index iff its slope has bounded partial coefficients. Berstel [B99] verified [M89] with some additional properties on the standard words. Carpi and de Luca [CD00] showed an explicit formula on the partial coefficients of the slope of a given Sturmian word to return the fractional index of the word and in process extended the same result on the Fibonacci words in [MP92]. This result in [CD00] is independently shown for integer indexes [V00] and verified with alternate proofs in [CW03, DL02]. [JP01] computed the maximum fractional power of runs that have a given finite factor as their period in a Sturmian word. Initial indexes— restricting the problem of maximal powers to the prefices of the words are studied in another line of research in this area. Allouche et. al [ADQZ01] and Damanik et. al [DKL00] (cf. [ADQZ01]) showed independently that each Sturmian word has an arbitrarily long prefix run. Berthé et. al [BHZ06] computed the initial fractional indexes of the runs in Sturmian words and thereby showed the following:
  - when the word is standard, the initial index is exactly one less than the index of the word itself,
  - a precise description of slopes $\theta$ admitting initial fractional indexes of 2.0 in their factor sets $T_\theta$.

In this paper, we characterize the runs in a Sturmian word and propose a linear time algorithm for their computation which, to the best of our knowledge, is an original study in literature. Following are the comparison of our results to those in literature:

- extends the solution of the same problem on Finonacci words in [IMS97] to the class of Sturmian words,
- generalizes the results in [DL03] to runs and linear-time computation of repetitions and their positions in a given Sturmian word,
- shows in process the number of *distinct* runs while this count is shown for the occurrences of runs in [BPR08] and for distinct powers of maximal and non-maximal repetitions in [DL03],
- verifies the results on integer indexes of Sturmian words [V00]. Also note the study of these results in context of fractional indexes in [CD00, CW03, DL02].

For a preliminary version of this paper, see [K98b, K99]. For a precise description of an algorithm to compute the runs in a given Sturmian word in linear-time as described in this paper, see [K99].

## 3. FINDINGS

In this section, the findings leading to Theorem 1 show how to detect the runs as they appear originally in their level. In consequence to Theorem 1, we show the cardinality of the set of distinct runs in a Sturmian word. Theorem 2 completes the study of runs in Sturmian words by the analysis of their extensions. Theorem 3 combines the results of the previous two theorems to present an alternate solution to the problem of integer indexes of Sturmian words.

**Lemma 1:**
Let $u$ be a binary word and $(p,p')$ any valid parameter pair. Then, $u=g^k$ if and only if $\alpha_{(p,p')}(u) = G^k$ where $k \geq 1$, $G = \alpha_{(p,p')}(g)$, $g$ and $G$ are primitive.

**Proof:** Let $u=g^k$ where $g$ is primitive, $k \geq 1$. Then, $\alpha_{(p,p')}(u) = \alpha_{(p,p')}(g^k) = [\alpha_{(p,p')}(g)]^k = G^k$ where $G = \alpha_{(p,p')}(g)$. $G$ is necessarily block-complete, thus $G$ can not be a repetition without $g = \alpha_{(p,p')}^{-1}(G)$ being a repetition, concluding that $G$ is primitive.

Conversely, let $\alpha_{(p,p')}(u) = G^k$ where $G$ is primitive and $G^k$ is necessarily block-complete. Then, $G$ has a full-block suffix and full-block prefix, and hence is a concatenation of full-blocks in all $k$ of its instances in $G^k$. Then, $\alpha^{-1}$ is defined on $G$ so that

$$u = \alpha_{(p,p')}^{-1}(G^k) = [\alpha_{(p,p')}^{-1}(G)]^k = g^k.$$

This implies $G = \alpha_{(p,p')}(g)$ and $g$ being a repetition would result $G$ being a repetition as well. □

We say that a run $u=(g,k,e)$ in a Sturmian word $S_i$ is a *reflection* of a run $u'=(g',k',e')$ in a level $j<i$ if an instance of $g$ is a left or a right shift of an instance of $\alpha_i \cdot \alpha_{i-1} \cdot \ldots \cdot \alpha_{j+1}(g')$. We say that a run $u$ is *original* in level $i$ if $u$ is not a reflection of a run of a previous level $j<i$. With this, whenever two or more instances of a period $g$ of a left or a right shift of a run $(g,k,e)$ is block-complete in level $i$, then the run $u$ is a reflection.

**Lemma 2:**
Let $u=g^k$ with $(g,k,e)$ be a run original in $S_i$ where $k>2$. Then, $g=a$.

**Proof:** Suppose $g>a$. Then $g$ contains at least one letter $b$ in order to be primitive. Let $g=vF_{(h)}w$ where $v$ and $w$ each equivalent to possibly empty partial blocks in the first instance of the period and $F_{(h)}$ a concatenation of $h$ full blocks, $h \geq 0$.

Since $g^2 = vF_{(h)}wvF_{(h)}w$ is a factor of $S_i$ and $S_i$ is a concatenation of full blocks, $wv$ is equal to a full block (short or long) when at least one of $v$ and $w$ is non-empty word, and $wv=\varepsilon$ when $v=w=\varepsilon$. Then, ($|v|$)-right shift of $u$ results in the run $u'=(F_{(h)}wv)^{k'}$ with $(F_{(h)}wv, k', e'_1, e'_2)$ where $e'_1=v$, $e'_2=F_{(h)}we$ if $|e|<|v|$, $e'_2=e[|v|+1..|e|]$ if $|e|\geq|v|$, $2 \leq k-1 \leq k' \leq k$ and the period $g'$ is block complete in all instances in the run. This implies that $u$ is a reflection of a run original at a previous level. $\square$

At this point we conclude that a run of order $k>2$ and period $g>a$ (i.e. $g \in \{a,b\}^* \setminus \{a, \varepsilon\}$) in $S_i$ is a reflection of a run that is original in an $S_j$ for some $j<i$.

**Eg.:**
a.) Consider the run $(ab,3,a)$ in Sturmian word $S_i=aababababaa$. Looking into $S_i$, $(p_i,p'_i)=(1,2)$, i.e., $\alpha_i(a)=ab$ and $\alpha_i(a)=aab$ since $ab$ is the repeating block. Furthermore, $S_i$ is necessarily suceeded by a letter $b$ implied by its suffix $aa$ and the existence of $bab$ sequence in the word. Thus, $S_ib$ is a full-block sequence in $S_i$. Looking into this run in level $i-1$,

$$\alpha_i^{-1}(S_ib) = baab$$

and the run $(ab,3,a)$ in level $i$ is the reflection of the original run $(a, 2, \varepsilon)$ in level $i-1$.

b.) Consider the run $(ab,2,a)$ in word $S_i=aababaa$. $S_i$ is the factor of any Sturmian word where $p'_i=1$. $S_i$ is also the factor of the set of Sturmian words where $p_i=1$ and either one of $p_{i-1}$ or $p'_{i-1}$ is 1. In $(ab,2,a)$, the period $ab$ is a full block in its second instance but a partial one in the first. Thus in $S_{i-1} = \alpha^{-1}_{(p,p')}(S_i)$, only the second instance of the period exists (as letter $a$ if $p_i=1$, letter $b$ otherwise) and $(ab,2,a)$ is an original run in $S_i$. $\square$

**Lemma 3:**
Let $u=g^2$ be a repetition original in $S_i$ where $g>a$. Then, $g=vT^hw$ and is a factor of $T^m wvT^n$ where $wv,T \in \{\alpha_i(a), \alpha_i(b)\}$, $wv \neq T$, and $h \geq 0$.

**Proof:**
Similar to the case in the previous proof, $g>a$ implies that $g=vF_{(h)}w$ where $v$ and $w$ are partial blocks and $F_{(h)}$ is a sequence of $h \geq 0$ full blocks in $S_i$.

Suppose first that $wv=\varepsilon$. Then, $g$ is block complete in both of its instances and the repetition $u$ is a reflection of a repetition that is original at an earlier level $S_j$ where $j<i$. Therefore, $wv > \varepsilon$. Since $g^2 = vF_{(h)}wvF_{(h)}w$ is a factor of $S_i$ and $S_i$ is a concatenation of full blocks, the instance of $X=wv$ in $g^2$ is a full block.

Let $S_i = s_1cv\ F_{(h)}X\ F_{(h)}wd\ s_2$ where $cv$ and $wd$ are full blocks. Then, $v \neq d$, otherwise ($|v|$)-right shift of $u$ resulted in a repetition of which the period is block complete in both of its instances and hence $u$ would be a reflection. Similarly, $c \neq w$. This implies that $cv \neq X$ and $wd \neq X$ and $cv=wd=X'$ where $(X,X')=(\alpha_i(l), \alpha_i(l'))$, $l,l' \in \{a,b\}$, $l \neq l'$.

Looking at $S_i = s_1X'F_{(h)}XF_{(h)}X'\ s_2$, we see that $F_1 = F_h$ since both full blocks $X$ and $X'$ cannot repeat in $S_i$. Let $Y = F_1 = F_h$ where $F_i$ denotes the i-th full block in $F_{(h)}$. Suppose $F_{(h)} \neq Y^h$ and thus $Y'$ is a factor of $F_{(h)}$ where $(Y,Y')=(X,X')$ or $(Y,Y')=(X',X)$. Let $F_j$, $j=2..h-1$ be the first occurrence of $Y'$ in $F_{(h)}$. Then, we are looking at:

$S_i = s_1X'Y^{j-1}Y'F_{j+1,h-1}YXY^{j-1}Y'F_{j+1,h-1}YX'\ s_2$

where $F_{j+1,h-1}$ denotes the $(j+1)$-st through $(h-1)$-st full blocks in $F_{(h)}$.

Considering the two possibilities:

Case (a): $X=Y$: Then,
$S_i = s_1Y'Y^{j-1}Y'F_{j+1,h-1}Y^{j+1}Y'F_{j+1,h-1}YY'\ s_2$.

We see that the block $Y$ repeats both $j-1$ and at least $j+1$ times between occurrences of $Y'$ which is not possible.

Case (b): $X=Y'$: Then,
$S_i = s_1Y^jY'F_{j+1,h-1}YY'Y^{j-1}Y'F_{j+1,h-1}Y^2\ s_2$.

In this case, $Y$ is the repeating block in $S_i$ since $j \geq 2$. Furthermore, $Y$ repeats either $j$ or $j-1$ times between every two occurrences of $Y'$ since both $Y^j$ and $Y'Y^{j-1}Y'$ are factors of $S_i$. From this, the words $Y$ and $Y^2$ following the two instances of $Y'\ F_{j+1,h-1}$ imply that $F_{j+1,h-1}=KY^{j-2}$ where $K=\varepsilon$ or $K$ is block complete in $S_i$ and $Y'$ is a proper suffix of $K$. This results in the sequence $Y'Y^{j-1}Y'Y^{j-1}Y'$ in $S_i$ and hence $\alpha_{i-1} \bullet \alpha_i(a) = Y^{j-1}Y'$ but not $Y^jY'$ repeats in $S_i$. From this and since $Y^jY'K$ is a sequence in $S_i$, $F_{j+1,h-1}$ has the prefix $Y^{j-1}Y'$, or, again, $K=\varepsilon$ with $F_{j+1,h-1}=Y^{j-2}$. Then, we have

$S_i = s_1\ Y^jY'K\ Y^{j-1}Y'\ Y^{j-1}\ Y'\ K\ Y^j\ s_2$.

Looking at $S_{i-2}$

$S_{i-2}= \alpha^{-2}(S_i) = \alpha^{-2}(s_1)\ b\ \alpha^{-2}(K)\ aa\alpha^{-2}(K)\alpha^{-2}(Y^js_2)$.

We know that either $K=(Y^{j-1}Y')^k$ for some $k \geq 0$ or $K$ is a concatenation of words $Y^{j-1}Y'$ and $Y^jY'$. If $K=(Y^{j-1}Y')^k$, then the period $g=v(Y^{j-1}Y')^{k+1}Y^{j-1}w$ of the repetition $u$ turns out to be a repetition itself with period $vY^{j-1}w$ and hence $u$ is not a primitive repetition, a contradiction to the condition of the lemma.

Suppose now that $K$ is a concatenation of words $Y^{j-1}Y'$ and $Y^jY'$. Since $\alpha^{-2}(Y^jY')=b$, $Y^jY'$ is neither a prefix or suffix of $K$ otherwise there is a $bb$ sequence in $S_{i-2}$. Let $K=(Y^{j-1}Y')^{t_1} Y^jY'K'(Y^{j-1}Y')^{t_2}$ where $t_1, t_2 \geq 1$, and $K'=\varepsilon$ or $K'$ is a sequence of $Y^{j-1}Y'$ and $Y^jY'$ where $Y^jY'$ is a suffix of $K'$. Then,

$S_{i-2} = \alpha^{-2}(s_1)ba^{t_1} b\alpha^{-2}(K')a^{t_1+t_2+2} b\alpha^{-2}(K')a^{t_2} \alpha^{-2}(Y^j)\alpha^{-2}(s_2)$.

$S_{i-2}$ would be unbalanced with $t_1$ and $t_1+t_2+2$ occurrences of letter $a$. □

**Corollary 1:**
Let $u=g^2$ be a repetition original in $S_i$ where $g>a$. Then, $u$ is exactly in one of the following two forms:

a.) a factor of $\alpha_i(bab)$ where $\alpha_i(a) = wv$, and $v$ and $w$ are respectively a prefix and a suffix of $g$,

b.) a factor of $\alpha_i(a^m ba^n)$ for some $m,n \geq 1$ where $\alpha_i(b) = wv$, and $v$ and $w$ are repectively a prefix and a suffix of $g$.

**Proof:**
By Lemma 3, $u = vX^k X'X^k w$ in $S_i$ where $v$ and $w$ are partial blocks with at most one of $v$ and $w$ is the empty word, $wv=X'$ and $(X,X')=(\alpha_i(l), \alpha_i(l'))$ with $l,l' \in \{a,b\}$, $l \neq l'$, $k \geq 0$. This sequence of $u$ is only possible as a reflection of $a^m ba^n$ or $bab$ for $m,n \geq 1$ and $k \leq \min\{m,n\}$ in $S_i$ □

At this point, we are ready for the following result which detects the runs at their original level in a Sturmian word.

**Theorem 1:**
Let $u=(g, k, e_2)$ be a run original in a Sturmian word $S_i$ of level $i$. Then, $u$ is exactly in one of the following forms:

i.)  $g=a$, $e_2=\varepsilon$, order $k$ for every $ba^k b$ sequence in $S_i$ where $k=p_i$ or $k=p'_i$ and $k>1$,

ii.)  $g=a^{p'_i} b(a^{p_i} b)^j$ where $j=0..\min\{m-1,n-1\}$, $k=2$, $e_2= a^{p'_i}$ for every $ba^m ba^n b$ sequence in $S_{i-1}$ $(m,n \geq 1)$ iff $p_i > p'_i$, iii.)  $g=(a^{p_i} b)^j a$ where $j=1..\min\{m,n\}$, $k=2$, $e_2 = a^{p_i-1}$ for every $ba^m ba^n b$ sequence in $S_{i-1}$ $(m,n \geq 1)$ iff $p_i < p'_i$, iv.)  $g= a^{p_i} b$, $k=2$, $e_2 = a^{p_i}$ for every $bab$ sequence in $S_{i-1}$ iff $p_i < p'_i$, v.)  $g= a^{p'_i} ba$, $k=2$, $e_2 = a^{p'_i-1}$ for every $bab$ sequence in $S_{i-1}$ iff $p_i > p'_i$,

**Proof:**
By Lemma 2, a repetition of order 3 or greater in $S_i$ which has a period greater than $a$ is a reflection. Lemma 2 also supports case (i) of this theorem.

By Lemma 3 and Corollary 1, a run of order 2 original in $S_i$ is a factor of $X^m X' X^n$ where $(X, X') = (\alpha_i(a), \alpha_i(b))$ and $m, n \geq 1$, or $(X, X') = (\alpha_i(b), \alpha_i(a))$ and $m = n = 1$, in both cases a prefix and a suffix of $X'$ is respectively a suffix of the first instance and the prefix of the second instance of the period in $u$. Let $H$ and $L$ represent respectively the short and long blocks in level $i$.

There are two possibilities:

Case (a): $X = H$
Case (b): $X = L$

Case (a): $X = H$: Then, by Corollary 1, there exists a repetition $u = g^2$ with period $g = vH^j w$ where $u$ is a factor of $Q = X^m X' X^n = H^m L H^n$, $wv = L$ and $j \leq \min\{m, n\}$. From this, there are $\min\{m, n\}$ repetitions as such for all $j = 1..\min\{m, n\}$.

If $v = \varepsilon$, then $w = L$ and $g = H^j L$. This is possible when $u$ is overlapping but not a factor of $Q$ and $u$ is block complete – a contradiction. So, $v \neq \varepsilon$. With a similar argument, $v \neq L$. Thus, $w \neq \varepsilon$ and $w[|w|] = a$.

Let $cud = cvH^j wvH^j wd$ so that $cv$ and $wd$ are block-complete. Since $e_1 = \varepsilon$, $c[|c|] \neq g[|g|]$. We know that $g[|g|] = w[|w|] = a$. Therefore, $c[|c|] = b$, and $v$ is equal to a full block sequence which is necessarily a short block since $w \neq \varepsilon$. Thus, $v = a^{p_{\min_i}} b$, $w = a$ and $g = (a^{p_{\min_i}} b)^j a$. Note here that, $v = a^{p_{\min_i}} b$ does not imply that $u$ is a block complete repetition since $v$ is necessarily a partial block in the second occurrence of the period in $u$.

Since $wd$ is a short but not a long block otherwise $u$ is the $(p_{max})$-left shift of the repetition with period $g' = H^{j-1} L$, order 2 and left extension $a^{p_{\min_i}} b$ with period block-complete in both of its instances and $u$ being a reflection in $X$. With this the right-extension $u$ is $e_2 = a^{p_{\min_i} - 1}$ since $a^{p_{\min_i} - 1}$ is the longest common prefix of $g$ and $d$. Then we have the repetitions $(g, 2, e_2) = ((a^{p_{\min_i}} b)^j a, 2, a^{|H|-2})$ for all $j = 1..\min\{m, n\}$ original in $S_i$, confirming cases (iii) and (v) of theorem.

Case (b): $X = L$: Then again, there exists a repetition $u = g^2$ with period $g = vL^j w$ where $u$ is a factor of $Q = X^m X' X^n = L^m H L^n$, $wv = H$ and $j \leq \min\{m, n\}$.

If $v = \varepsilon$ then $w = H$, implying that $g = L^j H$ but this is not possible since $(L^j H)^2$ is not a factor of $Q$. So, $v \neq \varepsilon$ and consecutively $v[|v|] = b$.

Looking into $S_i = cud = cvL^j HL^j wd$, $c[|c|] \neq b$ otherwise the repetition $u$ is not original in $S_i$ since $v$ is a full short block and the period is block complete in both of its occurrences. Since $e_1 = \varepsilon$ by condition and $c[|c|] = a$, $g[|g|] = b$. From this, $w = \varepsilon$ since $w$ is not a full block. Thus, $wv = v = a^{p_{\min_i}} b$.

Note that $v$ having a short block sequence does not imply that $u$ is block complete since $v$ is a partial block in the first instance of $g$. Since the suffix of first long block in $L^j H L^j$ is in $u$, $g = a^{P_{min_i}} b (a^{P_{max_i}} b)^k$ is a repetition original in $S_i$ where $k<j$.

Looking into the right extension, $a^{P_{min_i}}$ is the longest common prefix of $g$ and the sequence following the repetition $u$ in $S_i$, so $e_2 = a^{P_{min_i}}$. Then we have the repetitions $(g, 2, e_2) = (a^{P_{min_i}} b (a^{P_{max_i}} b)^l, 2, a^{P_{min_i}})$ for all $0 = 1..\min\{m,n\}-1$ original in $S_i$, confirming cases (ii) and (iv) of theorem. □

The following proof is a consequence of this theorem.

**Corollary 2:**

Let $R_u = \{(g_i, k_i, e_i) | (g_i, k_i, e_i) \text{ is a run in word } u\}$ and $G_u = \{g_i | (g_i, k_i, e_i) \text{ is a run in word } u\}$ denote respectively the set of distinct runs and the set of distinct periods in a word $u$. Let also $R'_u = \{(g_i, k_i, e_i) | (g_i, k_i, e_i) \text{ is a run in word } u \text{ and } (g_i, k_i, e_i) \text{ has an occurrence as a proper infix in } u\}$.

Let $s = \alpha_\pi(w)$ where $\pi = \{(p_1, p'_1), (p_2, p'_2), \ldots, (p_n, p'_n)\}$, $|w| \geq 2$ and $|w|_b \geq 1$, $w$ (and thus $s$) is a finite Sturmian word and the number of runs in $w$ is bound by a constant. Then, the following are true:

a.) $|G_s| = |G_w| + n + \sum_{i=1}^{n-1} p_i$

b.) $|R'_s| = |G_s| + n - k = |R'_w| + 2n - k + \sum_{i=1}^{n} p_i$ where $k=1$ if $\min\{p_n, p'_n\} = 1$, $k=0$ otherwise.

c.) $|R_s|, |G_s| \in \theta(\sum_{i=1}^{n} p_i)$. Furthermore, $|R_s|, |G_s| \in \Omega(n|w|^{1/n})$ and $|R_s|, |G_s| \in O(n+|w|/2^n)$.

**Proof:**
Each level $i$ with parameter pair $(p_i, p'_i)$ originates $p_i$ distinct runs in the next level due to cases (ii) and (iii) of Theorem 1, and any pair of these runs have distinct periods. The rest of the cases of the Theorem 1 define exactly 2 additional runs each with period $a$ originating in level $i$ (1 if $i=n$ and $k=1$, $k$ as in case (b) above). These discussions justify cases (a) and (b) of this corollary.

By cases (a) and (b), $|R'_s|$ and $|G_s|$ are in $\theta(n + \sum_{i=1}^{n} p_i)$ and, since $p_i \geq 1$, in $\theta(\sum_{i=1}^{n} p_i)$. $\sum_{i=1}^{n} p_i$ is at the peak (minimum) when $p_i = p_j$ for all possible $i$ and $j$ in range $[1,n]$, justifying $\Omega(n\sqrt[n]{|s|})$. At the other extreme, $\sum_{i=1}^{n} p_i$ is maximum when exactly one of the $p_i$ values is maximized, justifying $O(n + \frac{|s|}{2^n})$. The set $R_s \setminus R'_s$ contains those runs which are prefixes and suffixes of $u$. Such runs of each level are at most as many as the original runs of that level and $|R_s \setminus R'_s| \leq |R'_s|$. Thus, $|R_s|$ and $|R'_s|$ are asymptotically equivalent and $|R_s| \in \theta(|R'_s|)$, implying $|R_s| \in \Omega(n\sqrt[n]{|s|})$ and $|R_s| \in O(n + \frac{|s|}{2^n})$.
□

**Eg.:**
a.) Consider the finite Sturmian word *X=aaabaabaabaaabaabaabaab*. Since *aab* is the repeating block with two consecutive repetitions between the other block, *X* is necessarily followed by *aaab* in its superword and α$^{-1}$(*Xaaab*) = *baabaaab*.

The following two repetitions are original in *X,* exemplifying Case (*iii*) of Theorem 1.

* (*aaba*, 2, ε, *a*),
* (*aabaaba*, 2, ε, *a*).

Both of these repetitions are "centered" around the 12th position of *X*, i.e., the second occurrence of their period start at *X*[12]. These repetitions are original in *X* since no left or right shift of them is block-complete in both of the occurrences of the period. These repetitions are detected from the *aabaaa* sequence in α$^{-1}$(*X*). α$^{-1}$(*b*)=*aaab* in this sequence is overlapped by both occurrences of the period. Note that the repetition with period *aabaabaaba* and order 2 "centered" around the same position of *X* is also in *X* but not original since it is a left shift of the repetition with a block-complete period in its both occurrences. Indeed, it is the reflection of the repetition with period *baa* in α$^{-1}$(*X*).

b.) Consider now the finite Sturmian word *X=baabaaabaaabaabaaabaaabaaab*, this time with the long block repeating in the current level of the word. Similar to the above case, *X* is necessarily followed by *aab* in its superword. Furthermore, letter *b* in the prefix of X is necessarily the suffix of a long block in its superword since *aab* is not the repeating block in *X*. Thus, α$^{-1}$(*aaaXaab*) = *abaabaaab*.

The following two repetitions are original in *X,* exemplifying Case (*ii*) of the Theorem.

* (*aab*, 2, ε, *aa*),
* (*aabaaab*, 2, ε, *aa*).

Likewise, both of these repetitions are "centered" around the the same position of *X* and are detected from the *aabaaa* sequence in the previous level of the word. Observe that in both of these repetitions the period has a sequence equivalent to a block complete sequence-- a short block in (*aab*, 2, *e*, *aa*) and a short and a long block in (*aabaaab*, 2, *e*, *aa*). However, in both of these cases, the *aab* sequence which is equivalent to a short block in the prefix of the period is the proper suffix of a long block in the first occurrence of the period and thus the period is not block complete in its first occurrence. On the other hand, the repetiton with period *aabaaabaaab* and order 2 "centered" at same position as the above two is not original in *X* since the short block prefix of the period is block-complete this time also in its first occurrence. □

Theorem 1 clearly states the order and extensions of the repetitions as they appear originally in a Sturmian word. The following proofs complete the findings to compute the order and extensions of reflections of repetitions in the ultimate level over the defining sequence of a Sturmian word.

**Lemma 4:**
Let $u=(g, k, e_1, e_2)$ be a repetition in a Sturmian word $S_i$, $i \geq 1$. Then, the exact reflection of this repetition in $S_{i+1}=\alpha_{i+1}(S_i)$ has left extension $a^k b \alpha_{i+1}(e_1)$ and right extension $\alpha_{i+1}(e_2) a^k$ where $k=\min\{p_{i+1}, p'_{i+1}\}$.

**Proof:**
Assume $S_i = v e_1 g^k e_2 w$ so that $v[|v|] \neq g[|g|-|e_1|]$ and $w[1] \neq g[|e_2|+1]$. The left and right extensions of $g$ in $S_{i-1}$ have respectively $\alpha_{i+1}(e_1)$ and $\alpha_{i+1}(e_2)$ sequences.

In addition, the unequal letters $v[|v|]$ and $g[|g|-|e_1|]$ in $S_i$ form long and short blocks in $S_{i+1}$ and extend the left extension by their longest common suffix which is $a^k$, $k=\min\{p_{i+1}, p'_{i+1}\}$. Likewise, the right extension gains the additional sequence $a^k$ which is the longest common prefix of the long and short blocks. □

**Corollary 3:**
Let $u=(g,k,e)$ be a run in a Sturmian word $S_i$, $i \geq 1$. Then, the order of reflection of $u$ in $S_{i+1}$ is
i.)    $k$ if $|g| > |e|+1$
ii.)   $k+1$ otherwise, i.e. if $|g| = |e|+1$.

**Proof:**
Since $u$ is a run, $|g| \geq |e|+1$ by definition. By Theorem 2, a repetition gains additional left and right extensions equuivalent to the sequence $a^k b a^k$ where $k=\min\{p_{i+1}, p'_{i+1}\}$. This yields an additional full block in level $i+1$. This block adds to the extension of the repetition to form a sequence equivalent to $\alpha_{i+1}(g)$ in $S_{i+1}$ only when $|g|-|e|=1$. □

Lemma 4 states a trivial property to identify the extensions of a repetition in its reflection in the next level. Theorem 2 below expands this finding for the extensions of a repetition in its reflection in any further level.

**Theorem 2:**
Let $u$ be a repetition with period $g$ and respectively the left and right extensions $e_1$ and $e_2$ in a Sturmian word $S_i$ where $i \geq 1$. Then, the reflection of $u$ in $S_{i+n}$ with period $\alpha_{(i)}(g)$ has left extension

$$a^{p_{\min_{i+n}}} b \alpha_{(n)}(a^{p_{\min_{i+(n-1)}}} b) \alpha_{(n-1)}(a^{p_{\min_{i+(n-2)}}} b) \ldots \alpha_{(2)}(a^{p_{\min_{i+1}}} b) \alpha_{(1)}(e_1)$$

and right extension

$$\alpha_{(1)}(e_2) \alpha_{(2)}(a^{p_{\min_{i+1}}}) \ldots \alpha_{(n)}(a^{p_{\min_{i+(n-1)}}}) a^{p_{\min_{i+n}}}$$

where

$$p_{\min_i} = \min\{p_i, p'_i\}.$$

and

$\alpha_{(k)} = \alpha_{i+n} \cdot \alpha_{i+(n-1)} \cdot \alpha_{i+(n-2)} \cdot \ldots \cdot \alpha_{i+k}$ for $k \leq n$.

**Proof:**
Follows from Lemma 4 that a repetition gains respectively the left and right extensions $a^{P_{\min j}} b$ and $a^{P_{\min j}}$ in addition to the image of extension at each further level $j$. □

Using the results of Theorem 2, the following theorem provides an alternative solution [CD00, V00, DL02, CW03] to the problem of integer index of a Sturmian word given its defining sequence.

**Theorem 3:**
Let $s$ be a finite Sturmian word characterized by the defining sequence $\pi = \{(p_1, p'_1), (p_2, p'_2), \ldots, (p_{n-1}, p'_{n-1}), (p_n, p'_n)\}$. Let $P = \max\{p_j, p'_j : j < n-1\}$, $P_{n-1} = \max\{p_{n-1}, p'_{n-1}\}$, and $P_n = \max\{p_n, p'_n\}$. Then, the integer index of $s$:

$k = \max \{P_n, P_{n-1}+1, P+2\}$

if there exists an instance of $p$ value pair $(p_i, p'_i)$ in $\pi$ where $P = \max\{p_i, p'_i\}$ and $p'_{i+1} > p_{i+1}$,

$k = \max \{P_n, P_{n-1}+1, P+1\}$

otherwise, i.e., if there only exist instances of $p$ value pairs $(p_i, p'_i)$ in $\pi$ where $P = \max\{p_i, p'_i\}$ and $p'_{i+1} < p_{i+1}$.

**Proof:**
By Corollary 3 and from Theorem 1, among the possible repetitions original in a given $S_i$, only the following reflect with order $k+1$ in $S_{i+1}$:
a.) $(a, k, \varepsilon)$ for every $ba^k b$ sequence in $S_i$ where $k = p_i$ or $k = p'_i$ and $k > 1$ (case (i) of Theorem 1)
b.) $(a^{p'_i} b, 2, a^{p'_i})$ for every *aba* sequence in $S_{i-1}$ when $p_i > p'_i$ (case (ii) of Theorem 1)
c.) $(a^{p_i} b, 2, a^{p_i})$ for a *bab* sequence in $S_{i-1}$ when $p_i < p'_i$ (case (iv) of Theorem 1).

In (b) and (c) above, $\alpha_{i+1}(g) > 1 + \alpha_{i+1}(e_2)$ for all possible cases and thus these repetitions do not gain additional periods in levels greater than $i+1$. For case (a), since the period itself is a full block, namely $\alpha_{i+1}(a)$ in $S_{i+1}$, there is "potential" to gain an additional period in $S_{i+2}$.

Looking into the possibilities for this case, a repetition $(a,k,\varepsilon)$ reflects in one of the following two forms in level $S_{i+1}$:
(i) $(a^{p_{i+1}} b, k+1, a^{p_{i+1}})$ when $p'_{i+1} > p_{i+1}$,
(ii) $(a^{p'_{i+1}} ba, k+1, a^{p'_{i+1}} -1)$ when $p'_{i+1} < p_{i+1}$.

Among these two cases, the $|g|-|e|=1$ is satisfied only in the former and the repetition $(a,k,\varepsilon)$

original in $S_i$ reflects in $S_{i+2}$ as a repetition of order $k+2$. The period is larger in length than a block and thus $|g|-|e|=1$ is not satisfied for any further level $j>i+2$ and this repetition reflects with order $k+2$ at all levels greater than $i+2$.

From this, among the $a^k$ sequences in $S_i$, $i>2$, those succeeded by $p_{i+1}<p'_{i+1}$ in $\pi$ reflect with order $k+2$ and the rest by $k+1$, concluding the theorem. □

According to this, a Sturmian word defined by the uniform sequence
$\pi = \{(p,p'), (p,p'), (p,p'), ... \}$

generates a Sturmian word with the integer index value of $(\max\{p,p'\}+1)$ if $p>p'$, and of $(\max\{p,p'\}+2)$ if $p<p'$.

The "minimal" Sturmian word with the lowest possible "density" of '$a$'s is over the defining sequence

$\pi = \{(1,2), (1,2), (1,2), ... \}$.

The integer index of this word is 4 since it includes the factor $(ababa)^4$ in $\alpha^2(baab)$ and does not include a repetition $u^5$ for a nonempty factor $u$. The Fibonacci word over the defining sequence $\pi = \{(2,1), (2,1), (2,1), ... \}$ retains a unique property by having the minimum integer index of value 3 accross all Sturmian words.

## 4. CONCLUSION

Consider the η mapping associated with a finite defining sequence $\pi=\{(p_1, p'_1), (p_2, p'_2), (p_3, p'_3), ..., (p_n, p'_n)\}$

$\eta_\pi(p_i) = (p_{i+1}+1)p_i$
$\eta_\pi(p'_i) = (p_{i+1}+1)p'_i$
$\eta_\pi(1) = p'_{i+1}+1$

applied exactly once to each $p_i$, $p'_i$ and 1 term for $i=1$ to $n-1$ starting with an initial expression which represents each sequence $a^k$ by $k$ and letter $b$ by 1 in summation (Table 1).

| Level – $i$ | Word – $S_i$ | η-expression of $S_i$ |
|---|---|---|
| - | $a^k b$ | $k+1$ |
| 1 | $(a^{p_1}b)^k a^{p'_1}b$ | $(p_1+1)k+p'_1+1$ |
| 2 | $[(a^{p_2}b)^{p_1} a^{p'_2}b]^k (a^{p_2}b)^{p'_1} a^{p'_2}b$ | $[(p_2+1)p_1+p'_2+1]k+(p_2+1)p'_1+p'_2+1$ |
| ... | ... | ... |

**Table 1.** The representation in η-expression of initial levels of Sturmian word defined over $p=\{(p_1, p'_1), (p_2, p'_2), ..., (p_n, p'_n)\}$.

Consider the arithmetic value of this expression. Each $p_i$ (respectively $p'_i$) magnitude corresponds to a $a^{p_i}$ (respectively $a^{p'_i} b$) sequence, each 1 magnitude to a letter $b$, and each $p_i$ (respectively $p'_i$) multiplier of a parantheses pair corresponds to the superscript of the word representated by the expression enclosed by that pair. According to this, in a given word $s= \alpha_\pi(a)$, $s[j]=b$ if and only if the $j$-th magnitude in $\eta(s)$ is represented by a 1. Thus, the $\eta$-expression represents the associated Sturmian word $s$ and its arithmetic value is the length of $s$.

By Theorem 1, each $b a^{p_i} b$ and $b a^{p'_i} b$ originates one repetition in level $i$ (in level $i+1$ when the sequence is $bab$), and each $ba^m ba^n b$ where $m,n \in \{p_i, p'_i\}$ originates $\min\{m,n\}$ repetitions in level $i+1$. We take the view of tracing the repetitions by the $p_i$ and $p'_i$ terms appearing in $\eta$-expression so that a $P$ term where $P$ is $p_i$ or $p'_i$ term indicates the following two sets of original repetitions;

a.) one repetition of the form (*i*) of Theorem 1. This repetition is original in level *i*.
( (*iv*) or (*v*) of Theorem if the sequence is *bab*. In this case the repetition is original in level *i*+1. )

b.) $p_i$ or $p'_i$ repetitions of the form (*ii*) or (*iii*) of Theorem 1, associated to the letter $b$ succeeding that $p_i$ or $p'_i$ term. This set of repetition(s) are original in level $i+1$.

For each appeareance of a $P$ term in the $\eta$-expression where the term is not multiplying a parantheses pair, we detect exactly 1 original repetition. Each such repetition occurs as many times as the multiplication of the coefficients of the parantheses pairs enclosing it. A $p_i$ or $p'_i$ term that is the coefficient of a paranthesis represents maximum of $1+p_i$ repetitions each occurring as many times as that term is multiplied similarly to the above case. Each paranthesis pair has eactly one inner 1 term. In the expansion of the $\eta$-expression, this 1 is the multiplier of the occurrence-count of the repetitions represented by the $P$ term of that paranthesis-pair. And, since the length of a word is the arithmetic value of its $\eta$-expression, the number of runs in a Sturmian word is linear in the word length—verfying [BPR08]. See [K99] for a linear-time algorithm to compute the runs in a Sturmian word.